\documentclass[12pt]{article}
\usepackage{amsmath,amssymb,amscd,amsthm,latexsym,amsfonts, epsfig, graphicx}

\usepackage{color}

\definecolor{gr}{rgb}   {0.,   0.8,   0. } 
\definecolor{bl}{rgb}   {0.,   0.5,   1. } 
\definecolor{mg}{rgb}   {0.7,  0.,    0.7}

\begin{document}

\numberwithin{equation}{section} 

\newtheorem{theorem}{Theorem}[section]
\newtheorem{question}{Question}
\newtheorem{conjecture}[theorem]{Conjecture} 
\newtheorem{lemma}[theorem]{Lemma}
\newtheorem*{claim}{Claim}
\newtheorem{corollary}[theorem]{Corollary}
\newtheorem{proposition}[theorem]{Proposition}

\theoremstyle{remark}
\newtheorem*{remark}{Remark}
\newtheorem*{remarks}{Remarks}
\newtheorem*{merci}{Acknowledgements}
\newtheorem*{defi}{Definition}

\newcommand{\dv}{\operatorname{div}}
\newcommand{\R}{\operatorname{Re}}
\newcommand{\supp}{\operatorname{supp}}
\newcommand{\dist}{\operatorname{dist}}
\newcommand{\Lip}{\operatorname{Lip}}
\newcommand{\Dom}{\operatorname{Dom}}
\newcommand{\diam}{\operatorname{diam}}
\newcommand{\epi}{\operatorname{Epi}}
\newcommand{\aver}[1]{-\hskip-0.46cm\int_{#1}}
\newcommand{\avert}[1]{-\hskip-0.38cm\int_{#1}}
\newcommand{\ind}{1\hspace{-2.5 mm}{1}}
\newcommand{\sign}{\operatorname{sign}}

\newcommand{\NN}{\mathbb{N}}
\newcommand{\RR}{\mathbb{R}}
\newcommand{\CC}{\mathbb{C}}
\newcommand{\mH}{\mathcal{H}}

\def\barint_#1{\mathchoice
            {\mathop{\vrule width 6pt
height 3 pt depth -2.5pt
                    \kern -8.8pt
\intop}\nolimits_{#1}}%
            {\mathop{\vrule width 5pt height
3 pt depth -2.6pt
                    \kern -6.5pt
\intop}\nolimits_{#1}}%
            {\mathop{\vrule width 5pt height
3 pt depth -2.6pt
                    \kern -6pt
\intop}\nolimits_{#1}}%
            {\mathop{\vrule width 5pt height
3 pt depth -2.6pt
          \kern -6pt \intop}\nolimits_{#1}}}


\title{Sobolev spaces on multiple cones}
\author{P. Auscher\footnote{Univ Paris-Sud,
Laboratoire de Math\'ematiques, UMR 8628, Orsay,  F-91405;}
\footnote{CNRS, Orsay,  F-91405}
\\\texttt{{\footnotesize pascal.auscher@math.u-psud.fr}} \and
N. Badr\footnote{
Universit\'e de Lyon; CNRS; Universit\'e Lyon 1, Institut Camille Jordan, 43 boulevard du 11
Novem\-bre 1918, F-69622 Villeurbanne Cedex, France.}
\\\texttt{{\footnotesize badr@math.univ-lyon1.fr}}
\\}
\date{revised, may 2010}
\maketitle

\abstract{The purpose of this note is to discuss how various Sobolev spaces defined on multiple cones  behave with respect to density of smooth functions, interpolation and extension/restriction to/from $\RR^n$.    The analysis  interestingly combines use of Poincar\'e inequalities and  of some Hardy type inequalities.  



\section{Introduction}

The theory of Sobolev spaces on domains of the Euclidean spaces is well developed and numerous works and books are available. For multi-connected open sets, there is apparently nothing to say. However, depending on the topology of the boundary, the closure of  the space of test functions (ie compactly supported in $\RR^n$) might be a subtle thing. We propose here to  investigate 
the Sobolev spaces on multiple cones with common vertex as unique common point of their boundaries. Surprisingly, we did not find a treatment in the literature. 

Our motivation comes from Badr's PhD thesis  where  interpolation results for Sobolev spaces on complete metric-measure spaces are proved  upon the doubling property  and  families of Poincar\'e inequalities.  A question remained unsettled, namely whether the result is sharp, that is whether the conclusion is best possible given the hypotheses.    Multiple (closed) cones are sets where doubling (for Lebesgue measure) holds and $L^p$-Poincar\'e inequalities hold for some but not all $p$, more precisely for $p$  greater than dimension.   Cones are therefore simple but important examples  for this matter. 
The study of Sobolev spaces  on such sets provides us with the positive answer to our question and, in addition, we complete the interpolation result in this specific situation. As we shall see, these Sobolev spaces can be identified with the closure of test functions in the classical Sobolev space on open multiple cones.

 For simplicity, we work on  $\Omega$ the Euclidean (double) cone defined by 
$x_1^2 + \ldots + x_{n-1}^2 < x_n^2$, $n\ge 2$, but all the material extends right away to multiple cones with common vertex point (see Section 7 for natural extensions), and the cones need not be of revolution type.   Consider $W^1_p(\Omega)$ the usual first order Sobolev space on $\Omega$ and $\widetilde W^1_{p}(\Omega)$, the closure of smooth compactly supported functions in $\RR^n$ in $W^1_p(\Omega)$ if $p<\infty$, and   the space of bounded and Lipschitz functions on $\Omega$ that extend continuously at the origin if $p=\infty$.

The question we ask is: how do  they behave with respect to density of smooth functions, interpolation and extension/restriction to/from $\RR^n$? 

Our results (Sections 2,3,4 and 5) exhibit the specific role of  the vertex point. This role translates into a critical exponent (equal to dimension) and the $L^p$ Sobolev spaces   have different behaviors with respect to the various actions listed above. The following list   illustrates their properties : \begin{enumerate}
      \item The space $\widetilde W^1_{p}(\Omega)$ coincides with $W^1_p(\Omega)$ if $1\le p\le n$ but is of codimension 1 in $W^1_p(\Omega)$  for $n<p\le \infty$. (Section 2)
  \item The  spaces $W^1_p(\Omega)$, $1\le p \le \infty$, form a real interpolation family  (Section 3)
  \item The  spaces $\widetilde W^1_p(\Omega)$, $1\le p \le \infty$, do not form a real interpolation family. To obtain such a family, one needs to replace $\widetilde W^1_{n}(\Omega)$ by  a strict and dense subspace of it described in the text. (Sections 3, 4)
  \item For $p\in [1,\infty]$, $p\ne n$, the restriction operator to $\Omega$ maps  $W^1_{p}(\RR^n)$ continuously onto $\widetilde W^1_{p}(\Omega)$ and there exists a common linear continuous extension operator from $\widetilde W^1_{p}(\Omega)$ to $W^1_{p}(\RR^n)$. For $p=n$, these results hold with   $\widetilde W^1_{n}(\Omega)$ replaced by the  strict and dense subspace mentioned above. In particular,  $\Omega$ has the extension property for $W^1_p$ if and only if $1\le p<n$. (Section 5)
\end{enumerate}

Of course, some of these results are known and we give references along the way when we have been able to locate them. But some results, like the interpolation results, are new. We also point out that we give two proofs of the interpolation result. Although the one presented in Section 5 using restriction/extension looks more natural to users of Sobolev spaces on subsets of the Euclidean space, we prefer the one done in Sections 3 and 4, because it is more in the spirit of analysis on metric spaces and contains ideas that  we believe could be used in this context elsewhere. In particular, a special feature is that it allows to pass below the Poincar\'e exponent threshold by using  Hardy type inequalities.

We shall make use of the Sobolev space   $H^1_{p}(X)$ arising from geometric measure theory on $X=\overline\Omega$. It is defined as the completion for the $W^1_p(\Omega)$ norm of the space of Lipschitz functions with compact support for $p<\infty$ and as the space of bounded and Lipschitz functions in $X$ for $p=\infty$. It is easy to show it agrees with $\widetilde W^1_{p}(\Omega)$ and it turns out that it will be easier to work with the former in Section 3.  Finally, we make a connection with the Hajlasz-Sobolev space $M^{1,p}(X)$. In particular, we will show (Section 5) that the Hajlasz-Sobolev space $M^{1,n}(X)$ is a strict subspace of $H^{1}_{n}(X)$, which can be surprising.

In Section 6, we shortly describe the situation pertaining to these questions for homogeneous Sobolev spaces.

\section{Density}

Let $1\le p\le \infty$ and $O$ an open set of $\RR^n$. Define $W^1_p(O)$ as the space of functions\footnote{We consider real functions but everything is valid for complex functions.} $f \in L^{p}(O)$ such that 
$$
\|f\|_{W^1_p(O)}= \|f\|_{L^p(O)} +\|\nabla f\|_{L^p(O)}  <\infty.
$$
The gradient is defined  in the distributional sense in $O$. For $p<\infty$, denote by $\widetilde W^1_p(O)$ the closure of the space of $C^\infty_0(\RR^n)$ (the subscript $0$ means compact support) functions restricted to $O$ in $W^1_p(O)$. Among classical texts, we quote 
\cite{adams, EG, M1, M2, stein, ziemer}.

If $n<p<\infty$, recall that the Morrey-Sobolev embedding implies that if $f\in W^1_p(\Omega)$ then $f$ is H\"older continuous on each connected component $\Omega_\pm$ of $\Omega$,  the half-cones defined by $x\in \Omega$ and $\sign(x_n)=\pm 1$. Hence $f$ has limits in $0$ from $\Omega_+$ 
and $\Omega_-$. These limits, which we call $ f(0^+)$ and $f(0^-)$, may be different. 

 This lemma is classical and we include a proof for convenience. 

\begin{lemma}\label{lemma:density} Let $1\le p<\infty$ and $f\in W^1_p(\Omega)$. Assume $1\le p \le n$ or $n<p<\infty$ and $ f(0^+)=f(0^-)=0$. Then there exists a sequence of $C^\infty_0(\RR^n)$ functions $(\varphi_k)$ with support away from $0$ such that 
$\|f-\varphi_k\|_{W^1_p(\Omega)}$ tends to $0$.
\end{lemma}

\proof 
First, it is enough to consider $f \in W^1_{p}(\Omega_{+})$, with $f(0)=0$ if $p>n$. 
Second, we may also assume $f$ bounded by using the truncations $f_{N}=h_{N}(f)$ and $N\to \infty$, where $h_{N}(t)=-N$ if $t\le -N$, $h_{N}(t)=t$ if $-N\le t\le N$ and $h_{N}(t)=N$ if $t\ge N$. 
Next, we claim that we can approximate $f$ by a function in $g\in W_{p}^1(\Omega_{+})$ supported away from a ball centered at $0$. 
Assuming this claim, it suffices to convolve this approximation with a smooth mollifying function which has compact support inside $\Omega_{+}$ to conclude.

It remains to prove the above claim. Take $\chi\in C^\infty_{0}(\mathbb{R}^n)$ a positive, radial function, bounded by 1,  supported in the unit ball with $\chi=1$ in the half-unit ball. For $\epsilon>0$,  define $\chi_{\epsilon}(x)=\chi(\frac{x}{\epsilon})$ and take $f_{\epsilon}(x)=f(x)(1-\chi_{\epsilon}(x))$. Every $f_\epsilon=0$ on the ball of radius $\epsilon/2$. We distinguish between 3 cases:

\paragraph{Case $1\leq p<n$:}  By dominated convergence $f_{\epsilon}$ converges to $f$ in $ L^{p}(\Omega_{+})$. For the gradient, as  $f$ is bounded and $\|\nabla \chi_{\epsilon}\|_{p} \le C \epsilon^{n/p - 1}$, we conclude that $\nabla f_{\epsilon}$ converges to $\nabla f$ in $ L^{p}(\Omega_{+})$.

\paragraph{Case  $p=n$:} The function $f_{\epsilon}$ does not converge to $f$ in this case and we have to modify the construction. For $0<\delta <1$, we introduce the function $\eta_{\delta }(x)= \frac{|\ln \delta |}{|\ln |x||}$ if $|x|\le \delta $ and $\eta_{\delta }(x)=1$ if $|x|>\delta $. Take $f_{\epsilon,\delta }= f \eta_{\delta }(1-\chi_{\epsilon})=f_{\epsilon}\eta_{\delta}$ with $\delta^k =\epsilon$ and $k>0$ large. It is easy to show that $f_{\epsilon,\delta }$ converges to $f$ in $L^n(\Omega_{+})$ for $\epsilon \to 0$ and any $k$ fixed. For the gradient, using $|1-\chi_{\epsilon}| \le 1$, we have
$$
|\nabla (f - f_{\epsilon,\delta }) | \le |(1-\eta_{\delta })\nabla f|+ |\chi_{\epsilon}\nabla f| + |f \nabla \eta_{\delta }| + | f \eta_{\delta }\nabla \chi_{\epsilon}|.
$$
We observe that we assumed $f$ bounded. A computation shows that $\|  \eta_{\delta }\nabla \chi_{\epsilon}\|_{n}$ is bounded by $C/k$. So we pick and fix $k$ big enough. Next, $\|\nabla \eta_{\delta }\|_{n}$ goes to 0 as $\epsilon \to 0$ and the remaining  term $\|(1-\eta_{\delta })\nabla f\|_{n}+ \|\chi_{\epsilon}\nabla f\|_{n} \to 0$ by dominated convergence.

\paragraph{Case $p>n$:} 
By dominated convergence, $f_\epsilon$ converges to $f$ in $L^{p}(\Omega_+)$. For the gradient, we have $ (1-\chi_{\epsilon})\nabla f $ converges to $\nabla f$ in $L^p(\Omega_+)$ by  dominated convergence. It remains  to prove that  $\|f\nabla \chi_{\epsilon}\|_{L^{p}(\Omega_+)}$ tends to  $0$. By Morrey's theorem and recalling that $f(0)=0$ we have for every $x\in \Omega_{+},\, |x|<\epsilon$, \begin{equation} \label{eq:Morrey}
\left|\frac{f(x)}{\epsilon}\right|\leq C\left(\frac{|x|}{\epsilon}\right)^{1-n/p}\left(\barint_{\hspace{-5pt}\{|y|<2\epsilon\}\cap\Omega_{+}}|\nabla f|^{p}(y)dy\right)^{1/p}.
\end{equation}
 This implies
\begin{align*}\int_{\Omega_+}|f(x)&\nabla \chi_{\epsilon}(x)|^pdx \leq \int_{\{|x|\leq \epsilon\}\cap\Omega_+}\left|\frac{f(x)}{\epsilon}\right|^{p}dx 
\\
&\leq\int_{\{|x|\leq \epsilon\}\cap\Omega_+}\left(\frac{|x|}{\epsilon}\right)^{p-n}dx\ \frac{1}{\epsilon^{n}}\int_{\{|y|\leq 2\epsilon\}\cap\Omega_+}|\nabla f|^p(y)dy
\\
&\leq C\int_{\{|y|\leq 2\epsilon\}\cap\Omega_+}|\nabla f|^{p}(y)dy.
\end{align*}
We conclude noting that the last integral converges to $0$ when $\epsilon \rightarrow 0$ by the dominated convergence theorem. \qed

\begin{remark}
The density of functions in $W_{n}^{1}(\Omega_\pm)$ supported away from a ball centered at $0$ was also proved in \cite[Lemma 2.4]{costabel} for the special case of  dimension $n$ equals 2 (We are thankful to Monique Dauge for indicating this work). Their proof applies \textit{mutatis mutandis} for dimensions higher.
\end{remark}

\begin{corollary}\label{cor:wptilde}  Let $1\le p \le \infty$.   
If $p\leq n$, $\widetilde W^1_p(\Omega)=W^1_p(\Omega)$ and if $n<p$, 
 $\widetilde W^1_p(\Omega)=\{ f\in W^1_p(\Omega) \, ; \,  f(0^+)=f(0^-)\}$, and hence is of codimension 1 in $W^1_p(\Omega)$. 
 \end{corollary}
 
\proof
For $1\leq p\leq n$, the equality  follows immediately from Lemma \ref{lemma:density}. 
Assume now $n<p\le \infty$. Trivially $\widetilde{W}_{p}^{1}(\Omega)\subset \left\lbrace f\in W_p^1(\Omega); f(0^+)=f(0^-)\right\rbrace$. Conversely let $f\in W_p^1(\Omega),\, f(0^+)=f(0^-):=f(0)$. Then $g=f-f(0)\chi$, with $\chi\in C^{\infty}_{0}(\mathbb{R}^{n})$ supported in the unit ball with $\chi \equiv 1$ in a neighborhood of $0$ verifies $g(0^+)=g(0^-)=0$. By lemma \ref{lemma:density}  for $p<\infty$ and by definition for $p=\infty$, this yields $g\in \widetilde{W}_p^{1}(\Omega)$ and therefore $f=g+f(0)\chi$.
 \qed

\section{Real interpolation}

As far as $W^1_p(\Omega)$ is concerned, we have if $1\le p\leq\infty$ that 
$W^1_p(\Omega)= W^1_p(\Omega_+) \oplus W^1_p(\Omega_-)$ using restriction to  $\Omega_\pm$ and extension by 0 from $\Omega_\pm$ to $\Omega$. That is, if $f\in W^1_p(\Omega)$, we write
\begin{equation}\label{eq:dec}
f= {\bf 1}_{\Omega_+}f + {\bf 1}_{\Omega_-}f . 
 \end{equation}

Since $\Omega_\pm$ is a Lipschitz domain, it is known \cite{devore} that the family  of Sobolev spaces
$(W^1_p(\Omega_\pm))_{1\le p\le \infty}$ forms a scale of interpolation spaces for the real interpolation method. Hence the same is true for  $(W^1_p(\Omega))_{1\le p\leq \infty}$.

There is a second chain of spaces  appearing in the axiomatic theory of Sobolev spaces on a metric-measure space (\cite{hajlasz1}, \cite{hajlasz2}, \cite{heinonen1}). Let $X$ be the closure of $\Omega$. Then  $X$ equipped with Euclidean distance and Lebesgue measure, which we denote by $\lambda$, is a complete  metric-measure space. The balls are the restriction to $X$ of Euclidean balls centered in $X$. 
For $1\le p <\infty$, we denote by $H^1_p(X)$ the completion for the norm $W^1_p(\Omega)$ of  $\Lip_0(X)$, the space of Lipschitz functions in $X$ with compact support. For $p=\infty$, we set $H^1_\infty(X)=\Lip(X)\cap L^{\infty}(X)$. Identifying a Lipschitz function on $\Omega$ with its unique extension to $X$, $H_{\infty}^{1}(X)=\left\lbrace f\in W_{\infty}^1(\Omega); f(0^+)=f(0^-)\right\rbrace$.   There are also other Sobolev spaces of interest, like the Hajlasz spaces $M^{1,p}(X)$. We shall come back to this in Section 5.

    We recall the definitions of doubling property and Poincar\'e inequality:
\begin{defi}[Doubling property] Let $(E,d,\mu)$ be a metric-measure space. One says that $E$ satisfies the doubling property $(D)$ if there exists a constant $C<\infty$ such that for all $x\in E,\, r>0$ we have
\begin{equation*}\tag{$D$}
0<\mu(B(x,2r))\leq C \mu(B(x,r)).
\end{equation*}
\end{defi}

\begin{defi}[Poincar\'{e} Inequality]\label{PG} A  (complete)  metric-measure space $(E,d,\mu)$ admits a  $q$-Poincar\'{e} inequality  for some $1\leq q<\infty$, if there exists a positive constant $C<\infty$, such that for every continuous function $u$ and upper gradient $g$ of $u$, and for every ball $B$ of radius $r>0 $ the following inequality holds:
\begin{equation*} \tag{$P_{q}$}
\left(\barint_{\hspace{-5pt}B}|u-u_{B}|^{q}\, d\mu\right)^{\frac{1}{q}}\leq Cr \left(\barint_{\hspace{-5pt} B}g^{q}\, d\mu\right)^{\frac{1}{q}}. 
\end{equation*}

\end{defi}

There are weaker ways of defining the Poincar\'e inequalities but it amounts to this one when the space is  complete. See \cite{hajlasz2} for more on this and definition of upper gradients. 
On $X$, $|\nabla u|$ is an upper gradient of $u$. 

Let us recall Badr's theorem in this context (\cite{badr1}, Theorem 7.11).  On a metric-measure space there is a definition  $H^1_{p}(E)$ for $p\le\infty$ which, for $X=E$, is equivalent to the one given here.

\begin{theorem}[Badr]\label{th:Badr} Let $1\le q_{0}<\infty$. Assume $(E,d,\mu)$ is a complete metric-measure space with the doubling property and $q$-Poincar\'e inequalities  with  $q>q_{0}$. Then for   $q_{0}<p_{0}<p_1\le \infty$ and $ 1/ p = (1-\theta)/p_0 +   \theta/p_1$, 
\begin{equation}
(H^1_{p_0}(E), H^1_{p_1}(E))_{\theta,p}= H^1_p(E).
\end{equation}
\end{theorem}

The  space $(X,d, \lambda)$ has the doubling property and, as shown in \cite{hajlasz2} p.17,   it supports a $q$-Poincar\'e inequality  if and only if $n<q$. Thus, $(H^1_p(X))_{n< p\le \infty}$ is a scale of interpolation spaces for  the real interpolation method. As  observed and proved in \cite{badr}, Chapter 4 (see also Section 9 of \cite{badr1}), with arguments we reproduce here, $H^1_{p}(X)=W^1_{p}(\Omega)$ when $1\le p <n$, and this allowed her to identify  $H^1_p(X)$ as the interpolation space $(H^1_{p_0}(X), H^1_{p_1}(X))_{\theta,p} $ when $1\le p_0<p<p_1\le \infty$ and $ 1 /p = (1-\theta)/p_0 +   \theta/p_1$  with the restriction that either $n<p$ or $p_1<n$.

  The missing cases are somehow intriguing and for the sake of  curiosity  we provide a complete picture in the following result. More interestingly, we provide two proofs that cover all cases at once.  

\begin{theorem}\label{th:interpolationhp} If $1\le p_0<p<p_1\le \infty$ and $ 1/ p = (1-\theta)/p_0 +   \theta/p_1$, then 
\begin{equation}
(H^1_{p_0}(X), H^1_{p_1}(X))_{\theta,p}= \begin{cases} H^1_p(X), & {\rm if} \ p\ne n,\\
\widehat H^1_n(X) , & {\rm if} \ p= n.
\end{cases}
\end{equation}
\end{theorem}

We shall see that $\widehat H^1_n(X)$ is a strict subspace of $H^1_n(X)$. This implies in particular that Badr's interpolation result is sharp in the class of Sobolev spaces on metric-measure spaces: in this example,   the infimum of Poincar\'e exponents is also the smallest exponent  $p_0$ for which the   family $(H^1_p(X))_{p_0< p\le \infty}$ is a scale of interpolation spaces for  the real interpolation method. Hence, she could not get a better conclusion in general.  See section 5 for a further discussion on this.

The space $\widehat H^1_n(X)$ will incorporate a sort of Hardy inequality with respect to the vertex point.  To describe it, we need the following definition.

\begin{defi} For a function $f\colon X \to \RR$, we define its radial part $f_r$ and its anti-radial part $f_a$ as follows:  $f_r(x)$ is  the mean of $f$ on the sphere $S_{|x|}$ of radius $|x|$ restricted to $\Omega$ with respect to surface measure and $f_a(x)=f(x)-f_r(x)$. 
\end{defi}

The number  $f_r(x)$ depends only on the distance of $x$ to the origin, hence the terminology radial (even if $\Omega$ is not invariant by rotations). But note that both $f_r$ and $f_a$ depend on $\Omega$. Note that $f\mapsto f_r$ is a  contraction on $H^1_p(X)$. Denote by $r\colon \RR^n \to \RR, r(x)=|x|$.

\begin{defi} $\widehat H^1_n(X)= \{ f \in H^1_n(X)\, ; \, f_a/r \in L^n(X)\}$ with norm
$$
\|f\|_{\widehat H^1_n(X)} = \|f\|_{H^1_n(X)} +\|f_a/r\|_{L^n(X)}.
$$
\end{defi}

The following example shows that  $\widehat H^1_n(X)$ is a strict subspace of $H^1_n(X)$. Assume $n=2$ and $\beta>0$,  and consider the function $f$ on $X$, supported on $r\le 1/2$, $C^\infty$  away from 0, which is   $\sign(x_2) |\ln r|^{-\beta}$ for $r\le 1/4$. It is easy to check that  
$f\in H^1_2(X)$ for all $\beta>0$. Clearly, $f=f_a$ and $f/r \in L^2(X)$  if and only if $\beta>1/2$. Hence for $0<\beta\le 1/2$ we have $f\notin \widehat H^1_2(X)$.

Before we move on, the relation between $H^1_p(X)$ and $W^1_p(\Omega)$ is the following.

\begin{lemma}\label{lemma:hp} For $1\le p \le \infty$, $H^1_p(X)=\widetilde W^1_p(\Omega)$ with the same norm. 
\end{lemma}

\proof The equality at $p=\infty$ is obvious. Assume next that $p<\infty$. It is clear that $\widetilde W^1_p(\Omega) \subset H^1_p(X) \subset W^1_p(\Omega)$. Thanks to Corollary \ref{cor:wptilde}, we have our conclusion if $1\le p\le n$. Assume further $n<p$. Then functions in $\Lip_0(X)$ satisfy $f(0^+)=f(0^-)$. Since $f\mapsto f(0^\pm)$ are continuous on $W^1_p(\Omega)$, this passes to $H^1_p(X)$. Applying again Corollary \ref{cor:wptilde}, we deduce that $H_{p}^{1}(X)\subset \widetilde{W}_{p}^{1}(\Omega)$. \qed

\bigskip

To prove our theorem, we first  introduce the following spaces.
\begin{defi} 
For $1\le p \le  \infty$, set $\widetilde H^1_p(X)= \{f\in H^1_p(X)\, ; \, f/r \in L^p(X)\}$ with norm
$$
\|f\|_{\widetilde H^1_p(X)}=\|f\|_{H^1_p(X)} +\|f/r\|_{L^p(X)}= \|f\|_{W^1_p(\Omega)} + \|f/r\|_{L^p(\Omega)}.
$$
\end{defi}

\begin{lemma} 
For $1\le p \le  \infty$, $\widetilde H^1_p(X)$ is a Banach space which can be identified isometrically to  $\{f\in W^1_p(\Omega)\, ; \, f/r \in L^p(\Omega)\}.$ 
\end{lemma}

\proof There is nothing to prove if $1\le p\le n$ thanks to Corollary \ref{cor:wptilde} and Lemma \ref{lemma:hp}. Assume next $n<p\le \infty$.  Let $f\in \widetilde H^1_p(X)$, then the restriction of $f$ to $\Omega$ belongs to 
$\{f\in W^1_p(\Omega)\, ; \, f/r \in L^p(\Omega)\}$. Conversely if $f \in  W^1_p(\Omega)$ and $ f/r \in L^p(\Omega)$, then $f$ has a unique extension to a H\"older (Lipschitz if $p=\infty$) continuous function in both $\overline{\Omega_\pm}$. The condition $f/r \in L^p(\Omega)$ forces $f(0^+)=f(0^-)=0$. Hence this extension is in $ H^1_p(X)$ and thus in $\widetilde H^1_p(X)$. \qed
\bigskip

The next result is the main step.

\begin{theorem}\label{th:interpolationhptilde}
The family $(\widetilde H^1_p(X))_{1\le p\le \infty}$ is a scale of interpolation spaces for  the real interpolation method.
\end{theorem}

This result is proved in the next section. We continue with

\begin{proposition}\label{tildehc} If $1\le p<n$, $\widetilde H^1_p(X)= H^1_p(X)$. If $n<p\le \infty$, $\widetilde H^1_p(X)=\{f \in H^1_p(X)\, ;  \, f(0)=0\}$ and has codimension 1 in $H^1_p(X)$.
\end{proposition} 
Before we prove this proposition we need the following  Hardy type inequality  (we thank Michel Pierre for indicating a simple proof):
\begin{lemma}\label{hardy}
Let $1\leq p\leq \infty$ with $p\neq n$. Then there exists a constant $C=C(p,\Omega)$ such that
\begin{equation}\label{uc}
\int_{\Omega}\left|\frac{f}{r}\right|^{p}dx \leq C\int_{\Omega}|\nabla f|^{p}\,dx
\end{equation}
for every $f\in H_p^1(X)$ with, in addition,  $f(0)=0$ if $p>n$,  and \eqref{uc} is understood with  $L^\infty$ norms if $p=\infty$.
\end{lemma}

The example above shows that the lemma is false when $p=n$.

\proof  Assume first $1\leq p<n$. Take $f\in \Lip_{0}(X)$. We have
\begin{align*}
\int_{\Omega_{+}}\left|\frac{f}{r}\right|^{p}dx&=\int_{\Omega_{+}\cap S_1}\int_{0}^{\infty}r^{n-1-p}|f(r,\theta)|^{p}drd\sigma(\theta)
\\
&=\int_{\Omega_{+}\cap S_1}\left[\frac{1}{n-p}r^{n-p}|f(r,\theta)|^{p}\right]_{0}^{\infty}d\sigma(\theta)
\\
&\quad -\int_{\Omega_+\cap S_1}\int_{0}^{\infty}\frac{1}{n-p}r^{n-p}p|f|^{p-1} \,\sign f\,\frac{\partial{f}}{\partial{r}}drd\sigma(\theta) 
\\
&=-\frac{p}{n-p}\int_{\Omega_+\cap S_1}\int_{0}^{\infty}\left|\frac{f}{r}\right|^{p-1}\,\sign f\,\frac{\partial{f}}{\partial{r}}r^{n-1}drd\sigma(\theta)
\\
&\leq \frac{p}{n-p}\left(\int_{\Omega_+\cap S_1}\int_{0}^{\infty}\left|\frac{f}{r}\right|^{p}r^{n-1}drd\sigma(\theta)\right)^{\frac{p-1}{p}}
\\
&\quad \times \left(\int_{\Omega_+\cap S_1}\int_{0}^{\infty}\left|\frac{\partial{f}}{\partial{r}}\right|^{p}r^{n-1}drd\sigma(\theta)\right)^{\frac{1}{p}}.
\end{align*}
After simplification, we get (\ref{uc}) on $\Omega_+$. We do the same for the integral on $\Omega_{-}$ and therefore (\ref{uc}) holds for every $f\in \Lip_{0}(X)$. By density, (\ref{uc}) holds for every $f\in H_{p}^{1}(X)$.
\\

Assume next $n<p<\infty$.
Let $f\in \Lip_{0}(X)$ such that $f(0)=0$.  We denote $A=\int_{\Omega_{+}\cap\{|x|>\epsilon\}}\big|\frac{f}{r}\big|^{p}dx$, where $\epsilon>0$. By Morrey's theorem, we have for every $x\in \Omega$, $|f(x)|\leq C\|\,|\nabla f|\,\|_{p}|x|^{\alpha}$ with $\alpha=1-n/p$.
 Repeating the computation of (\ref{uc}) and since $f$ has a compact support, one obtains \begin{align*}
A&=\int_{\Omega_{+}\cap S_1}\int_{\epsilon}^{\infty}r^{n-1-p}|f(r,\theta)|^{p}drd\sigma(\theta)
\\
&=\int_{\Omega_{+}\cap S_1}\left[\frac{1}{n-p}r^{n-p}|f(r,\theta)|^{p}\right]_{\epsilon}^{\infty}d\sigma(\theta)
\\
&-\int_{\Omega_{+}\cap S_1}\int_{\epsilon}^{\infty}\frac{1}{n-p}r^{n-p}p|f|^{p-1} \,\sign f\,\frac{\partial{f}}{\partial{r}}d\sigma(\theta) dr
\\
&=\frac{\epsilon^{n-p}}{p-n}\int_{\Omega_{+}\cap S_1} |f(\epsilon,\theta)|^p d\sigma(\theta)
+\frac{p}{p-n}\int_{\Omega_+\cap S_1}\int_{\epsilon}^{\infty}\left|\frac{f}{r}\right|^{p-1}\sign f\frac{\partial{f}}{\partial{r}}r^{n-1}drd\sigma(\theta)
\\
&\leq C^{p} \|\nabla f\|_{p}^{p}
\\
&+ \frac{p}{p-n}\left(\int_{\Omega_+\cap S_1}\int_{\epsilon}^{\infty}\left|\frac{f}{r}\right|^{p}r^{n-1}drd\sigma(\theta)\right)^{\frac{p-1}{p}}\left(\int_{\Omega_+\cap S_1}\int_{\epsilon}^{\infty}\left|\frac{\partial{f}}{\partial{r}}\right|^{p}r^{n-1}drd\sigma(\theta)\right)^{\frac{1}{p}}.
\end{align*}
This yields
\begin{equation*}\label{AB14}
A\leq C^{p}\|\,|\nabla f|\,\|_{p}^{p}+\frac{p}{p-n}A^{\frac{p-1}{p}}\|\,|\nabla f|\,\|_{p}.
\end{equation*}
Plugging $$A^{\frac{p-1}{p}}\|\,|\nabla f|\,\|_{p}\leq \frac{\delta^{p'} A}{p'}+\frac{1}{p\delta^p}\|\,|\nabla f|\,\|_{p}^{p}$$ for every $\delta>0$,  with $p'=\frac{p}{p-1}$,
one obtains 
 $$ A(1-\frac{p\delta^{p'}}{(p-n)p'})\leq (C^p+\frac{1}{(p-n)\delta^p})\|\,|\nabla f|\, \|_{p}^{p}.
 $$
 Choosing $\delta$ small enough, we deduce that
  $$\int_{\Omega_{+}\cap\{|x|>\epsilon\}}\left|\frac{f}{r}\right|^{p}dx\leq C\int_{\Omega_+}|\nabla f|^{p}dx.
  $$
  We then let $\epsilon\rightarrow 0$.
We do the same for the integral on $\Omega_{-}$ and therefore (\ref{uc}) holds for every $f\in \Lip_{0}(X)$ such that $f(0)=0$. By density, (\ref{uc}) holds for every $f\in H_{p}^{1}(X)$ such that $f(0)=0$.
\\

When $p=\infty$,  (\ref{uc}) is a direct consequence of the definition of $H_{\infty}^1(X)$ and that $f(0)=0$ with the mean value theorem. \qed

\proof[Proof of Proposition \ref{tildehc}] When $1\leq p<n$, Lemma \ref{hardy} shows that $H_p^1(X)\subset \widetilde{H}_p^1(X)$ and the proposition follows. Now, when $p>n$, Lemma  \ref{hardy} yields $\left\lbrace f\in H_{p}^{1}(X); f(0)=0\right\rbrace\subset \widetilde{H}_{p}^{1}(X)$. 
Conversely if $f\in \widetilde{H}_p^1(X)$, by the continuity of $f$ at $0$ and the $L^p$ integrability of $f/r$ we easily see that $f(0)=0$.

It remains to prove that $\widetilde{H}_{p}^{1}(X)$ is of codimension $1$ in $H_p^1(X)$. This follows by writing $f\in H_{p}^{1}(X)$ as $f=f-f(0)\chi+f(0)\chi$, where $\chi\in C^{\infty}_{0}(\mathbb{R}^{n})$, $\supp \chi\subset B(0,1)$ and $\chi=1$ in a neighborhood of $0$, and using the above characterization of $\widetilde{H}_p^1(X)$. \qed

Although this is a simple description of $\widetilde H^1_p(X)$, the jump at $p=n$ does not allow us to use this result to conclude for Theorem \ref{th:interpolationhp}. We need to further analyze 
the radial  and antiradial parts of a function.

\begin{lemma}\label{lemma:ra} Let $1\le p \le \infty$.
\begin{enumerate}
\item For a function $f$ depending only on the distance to the origin, $f\in H^1_p(X) \Longleftrightarrow f\in W^1_p(\RR^n)$ with same norm up to a constant.
\item Assume $p\ne n$. For a function $f\colon X \to \RR$ and $f_a=f-f_r$, we have $f_a\in \widetilde H^1_p(X)
\Longleftrightarrow f_a\in  H^1_p(X)$ with comparable norms.
\end{enumerate}
\end{lemma}

\proof The first item is trivial. The constant is the ratio of the surface measure  of $\Omega$  inside the unit sphere divided by the surface measure of the unit sphere. 

As for the second item, it follows  from the previous proposition directly if $p<n$ and by observing that $f_a(0)=0$
if $p>n$. \qed 

\bigskip
Let us recall the following definition:
\begin{defi} 
Let $f$  be a measurable function on a measure space $(X,\mu)$. The decreasing rearrangement of $f$ is the function $f^{*}$ defined for every $t\geq 0$ by
$$
f^{*}(t)=\inf \left\lbrace\lambda :\, \mu (\left\lbrace x:\,|f(x)|>\lambda\right\rbrace)\leq
t\right\rbrace.
$$
The maximal decreasing rearrangement of
$f$ is the function $f^{**}$ defined for every $t>0$ by
$$
f^{**}(t)=\frac{1}{t}\int_{0}^{t}f^{*}(s) ds.
$$
\end{defi}
\begin{remark}\label{rearrangement}
 It is known that $(\mathcal{M}f)^{*}\sim f^{**}$ with $\mathcal{M}$ the Hardy-Littlewood maximal operator, $ \|f^{**}\|_p \sim \|f\|_p$ for all $ p>1$ (see \cite{stein1}, Chapter V, Lemma 3.21, p.191 and Theorem 3.21, p.201) 
and $\mu (\left\lbrace x:\, |f(x)|>f^{*}(t)\right\rbrace)\leq t$ for all $t>0$. We refer to \cite{bennett}, \cite{bergh}  for other properties of $f^{*}$ and $f^{**}$.
\end{remark}

We can now complete the proof of  Theorem \ref{th:interpolationhp}.

\proof Let us examine the case where neither $p_0, \,p_1$ is $n$. By the reiteration theorem, this reduces further to $p_0=1,\,p_1=\infty$.  Set  $F_p=(H^1_{1}(X), H^1_{\infty}(X))_{\theta,p}$ with $\theta=1-1/p$.

Let $f\in F_p$. Since $f\mapsto f_r$ is contracting on $H^1_q(X)$ for all $1\le q\le \infty$ and using Lemma \ref{lemma:ra}, one has that  $$K(f_r,t, W^1_1(\RR^n), W^1_\infty(\RR^n)) \le C K(f,t, H^1_1(X), H^1_\infty(X)).$$
$K$ is the $K$-functional of interpolation defined as in \cite{bennett}, \cite{bergh}. Hence $f_r\in W^1_p(\RR^n)$ by classical interpolation for the $W^1_p(\RR^n)$. Thus  $f_r \in H^1_p(X)$ by Lemma \ref{lemma:ra}. We also have by Lemma \ref{lemma:ra} again, 
$$K(f_a,t, \widetilde H^1_1(X), \widetilde H^1_\infty(X)) \le C K(f,t, H^1_1(X), H^1_\infty(X)).$$
Theorem \ref{th:interpolationhptilde} shows then that $f_a \in \widetilde H^1_p(X)$. 
We conclude that $f\in H^1_p(X)$  if $p\ne n$  and $f \in \widehat H^1_n(X)$ if $p=n$.

Reciprocally, let $f\in H^1_p(X)$  if $p\ne n$ and  $f\in \widehat H^1_n(X)$ if $p=n$. By Lemma \ref{lemma:ra}, whatever $p$ is, we have that $f_r\in W^1_p(\RR^n)$ and $f_a \in \widetilde H^1_p(X)$. By Theorem \ref{th:interpolationhptilde},  $f_a \in  (\widetilde H^1_1(X), \widetilde H^1_\infty(X))_{\theta,p}$ with $\theta=1-1/p$. Hence $f_a \in F_p$. For the radial part, for each $t>0$, one can find   a decomposition $f_r=g_t+h_t$ almost  minimizing for $K(f_r,t,W^1_1(\RR^n),  W^1_\infty(\RR^n))$ and one can assume both $g_t$ and $h_t$ are radial. Thus Lemma \ref{lemma:ra} implies that $g_t\in H^1_1(X)$ and $h_t \in H^1_\infty(X)$, hence $f_r \in F_p$.

It remains to study the case where  $p_0$ or $p_1$ is equal to $n$. Let us consider the case $p_1=n$ as the other one is similar. It is also enough to look at the result when $p_0=1$. As we know all interpolation spaces between $H^1_1(X)$ and $H^1_\infty(X)$, by the reiteration theorem, if $1<p<n$ and $\frac 1 p = {1-\theta} + \frac \theta n $ we have $(H^1_{1}(X), \widehat H^1_{n}(X))_{\theta,p}= H^1_p(X)$. 
Hence, we have
$$
H^1_p(X)=  (H^1_{1}(X), \widehat H^1_{n}(X))_{\theta,p} \subset (H^1_{1}(X), H^1_{n}(X))_{\theta,p} \subset H^1_p(X).
$$
The last inclusion is the easy part of the interpolation: we recall that
 $$
K(f,t, H_{1}^{1}(X), H_{n}^{1}(X))\geq K(f,t,L_{1},L_{n})+K(|\nabla f|,t,L_{1},L_{n})
$$ 
and that
$$
K(f,t,L_{1},L_{n})\sim \int_{0}^{t^{\alpha}}f^{*}(u)\,du+t\left(\int_{t^{\alpha}}^{\infty}f^{*}(u)^{n}\,du\right)^\frac{1}{n},
$$
where $\frac{1}{\alpha}=1-\frac{1}{n}$.
Integrating, using  the definition of the maximal rearrangement function and its properties in the  Remark  above, we deduce that
\begin{equation*}
\|f\|_{(H^1_{1}(X), H^1_{n}(X))_{\theta,p}}^{p}=\int_{0}^{\infty}\left(t^{-\theta} K(f,t, H_{1}^{1}(X), H_{n}^{1}(X))\right)^{p}\frac{dt}{t}
\geq C\|f\|_{H_p^1(X)}^p
\end{equation*}
and therefore
$
(H^1_{1}(X), H^1_{n}(X))_{\theta,p} \subset H^1_p(X)$.

This concludes the proof.
\qed

\begin{remark} The inclusion $\widehat H^1_n(X) \subset H^1_n(X)$ is dense.  This is due to the fact that as $H^1_n(X) = W^1_n(\Omega)$,  the space of restrictions to $X$ of smooth functions on $\RR^n$ with compact support in $\RR^n\setminus \{0\}$, a subspace of 
$\widehat H^1_n(X)$, is dense in $H^1_n(X)$.  Note that the last part of the argument shows that $(H^1_{p_{0}}(X), H^1_{p_{1}}(X))_{\theta,p} = H^1_p(X)$ for the appropriate $p$ whenever $p_{0}$ or $p_{1}$ is equal to $n$. In particular, when $p_{0}=n$ this furnishes an endpoint to Badr's result for the case of $E=X$.  
\end{remark}

\begin{remark}
As $X$ is symmetric with respect to $S:x\mapsto -x$, we can define $\widehat H^1_n(X)$ differently by doing an analysis with even and odd parts.  Define the even and odd parts $f_e$ and $f_o$ of a function $f\colon X\to \RR$ as $f_{e}=\frac 1 2 (f+f\circ S)$ and $f_{o}=\frac 1 2 (f-f\circ S)$. We have that 
$\widehat H^1_n(X)= \{ f \in H^1_n(X)\, ; \, f_o/r \in L^n(X)\}$.
 Let $f\in H^1_n(X)$.  
 Write  $f_o= (f_r)_o+ (f_a)_o$ and easily $(f_r)_o=0$. Hence
$$
f_a \in \widetilde H^1_n(X) \Longrightarrow (f_a)_o \in \widetilde H^1_n(X)  \Longrightarrow f_o \in \widetilde H^1_n(X).
$$
Next, write $f_a=(f_e)_a + (f_o)_a$.  We claim that $(f_e)_a/r \in L^n(X)$.  Hence, 
$$
f_o \in \widetilde H^1_n(X) \Longrightarrow (f_o)_a \in \widetilde H^1_n(X)  \Longrightarrow f_a \in \widetilde H^1_n(X).
$$
To see the claim, we observe that  the evenness of $f_e$ implies that $(f_e)_r(x)$ is also equal to the mean of $f_{e}$  on $\Omega_+ \cap S_{|x|}$.  Thus $(f_{e})_{a}$ is an even function with mean value 0 on each $\Omega_{\pm} \cap S_{|x|}$. Applying the classical (with gradient instead of upper gradient)  Poincar\'e inequalities ($P_{n}$) for spherical caps (that is geodesic balls) with respect to surface measure on $S_{|x|}$, we obtain
$$
\int_{\Omega_\pm\cap S_{|x|}}  {|(f_e)_a|^n} d\sigma(\theta)  \le C(n,\Omega_\pm) |x|^n \int_{\Omega_\pm \cap S_{|x|}}  {| \nabla_\theta (f_e)_a|^n} d\sigma(\theta)
$$
where $\nabla_\theta $ is the tangential gradient on $S_{|x|}$. Notice that $| \nabla_\theta (f_e)_a | \le |\nabla (f_e)_a |$ on $\Omega \cap S_{|x|}$. So adding the two inequalities,  multiplying by $r^{-n}=|x|^{-n}$ and integrating with respect to $dr$  we obtain
$$
\int_{\Omega} \frac {|(f_e)_a|^n}{r^n} dx  \le  C(n,\Omega) \int_{\Omega} {|\nabla (f_e)_a|^n}dx.
$$
Hence $(f_e)_a/r \in L^n(X)$ as claimed.

 \end{remark} 
 
\begin{remark} We have used Poincar\'e inequalities ($P_{n}$) for spherical caps  on spheres. Equipped with geodesic distance and surface measure, they are spaces of homogeneous type. Actually, ($P_{p}$) hold for  all $p\ge 1$ and all geodesic balls (including the sphere itself):
$$
\int_{B}|f(\theta) -m_{B}f|^p d\sigma(\theta) \le c(p,n) \diam(B)^p \int_{B} | \nabla_\theta f(\theta)|^p d\sigma(\theta).
$$
We have not been able to locate this result explicitly in the literature, neither can we say who proved it first. But it is not a new fact. It can be obtained from  \cite{LS} seeing them as submanifold in $\RR^n$. One can also apply results in  \cite[Theorem B.10]{Se}. One can also relate this to isoperimetry  and Sobolev inequalities (see, e.g., \cite[Chapter IV]{C}), especially if, unlike us, one is after best constants. A pedestrian approach to prove Poincar\'e inequalities for spherical caps (or more general Lipschitz subdomains of the sphere) is to pullback integrals  $\int_{B}|f(\theta) -a|^p d\sigma(\theta)$, $a$ constant, via a stereographic projection and use Poincar\'e inequalities on balls (or bounded Lipschitz domains) of $\RR^n$. This  easily works if $B$ is contained in a hemisphere by choosing an opposite pole. If this is not the case, cut $B$ in two equal parts along an equator and use the argument above for each part, using again  Poincar\'e inequalities  for bounded Lipschitz domains of $\RR^n$.
\end{remark}

\section{Proof of Theorem \ref{th:interpolationhptilde}}
For the proof of Theorem \ref{th:interpolationhptilde}, we need a Calder\'on-Zygmund decomposition as in \cite{badr1}. We incorporate here a further control to take care of the vertex point.
 
Let $1<p<\infty$ and $f\in \widetilde{H}_{p}^{1}(X)$. Identifying  $f$ to its restriction to $\Omega$, write $f=f|_{\Omega_{+}}+f|_{\Omega_{-}}=f_{+}+f_{-}$. We establish the following Calder\'on-Zygmund decomposition for $f_{+}$ and the same decomposition holds for $f_{-}$.

\begin{proposition}
[Calder\'{o}n-Zygmund lemma]\label{CZ}  
Let  $\alpha>0$. Then one can find a collection of balls $(B_{i+})_{i}$ of $\Omega_{+}$, functions $b_{i+}$ and a Lipschitz function $g_+$ such that the following properties hold:
\begin{equation}
f_+ = g_{+}+\sum_{i}b_{i+}  \quad {\rm on}\  \Omega_{+}\label{df}
\end{equation}
\begin{equation}
|g_+(x)| + \frac{|g_{+}(x)|}{|x|} +  |\nabla g_{+}(x)|\leq C\alpha\quad \lambda-a.e\; x\in \Omega_+ \label{eg}
\end{equation}
\begin{equation}
\supp b_{i+}\subset B_{i+}, \,\barint_{\hspace{-5pt}B_{i+}}\left(|b_{i+}|+\frac{|b_{i+}|}{|x|}+|\nabla b_{i+}|\right)dx\leq C\alpha\label{eb}
\end{equation}
\begin{equation}
\sum_{i}\lambda(B_{i_{+}})\leq C\alpha^{-p}\int_{\Omega_+} \left(|f_{+}|+\frac{|f_{+}|}{|x|}+|\nabla f_{+}|\right)^p dx
\label{eB}
\end{equation}
\begin{equation}
\sum_{i}\chi_{B_{i+}}\leq N \label{rb}.
\end{equation}
The constants $C$ and $N$ only depend on $p$ and on the constants in $(D)$ and $(P_{1})$ in $\Omega_{+}$.
\end{proposition}

A ball of $\Omega_{+}$ is the restriction to $\Omega_{+}$ of an open ball of $\RR^n$ having  center in $\Omega_{+}$. 

\proof To simplify the exposition, we  omit the index $+$  keeping it only for $\Omega_+$. For $x\in \mathbb{R}^{n}$, denote $r(x)=|x|$. Consider
$$U=\left\lbrace x \in \Omega_+ : \mathcal{M}_{\Omega_+}(|f|+\frac{|f|}{r}+|\nabla f|)(x)>\alpha \right\rbrace$$ with
$$
\displaystyle\mathcal{M}_{\Omega_+}f(x)= \sup_{B:\,x\in
B}\frac{1}{\lambda(B)}\int_{B}|f|dx
$$
where $B$ ranges over all balls of  $\Omega_+$. Recall that $\mathcal{M}_{\Omega_+}$ is of weak type $(1,1)$ and bounded on $L^{p}(\Omega_+,\lambda),\,1<p\leq\infty$.
If $U=\emptyset$, then set
$$
 g=f\;,\quad b_{i}=0 \, \text{ for all } i
$$
so that (\ref{eg}) is satisfied according to the Lebesgue differentiation theorem. Otherwise the maximal theorem gives us
\begin{equation}\label{mO}
	\lambda(U)
			 \leq C \alpha^{-p} \int_{\Omega_{+}} \left(| f|  + \frac{|f|}{r}+|\nabla f|\right)^p dx <\infty
\end{equation}
 In particular $U \neq \Omega_+$ as $\lambda(\Omega_+)=\infty$. Let $F$ be the complement of $U$ in $\Omega_{+}$. Since $U$ is an open set distinct of $\Omega_+$, we use  a Whitney decomposition of $U$ (\cite{coifman}): one can find
pairwise disjoint balls  $\underline{B_{i}}$  of $\Omega_{+}$  and  two constants $C_{2}>C_{1}>1$,  such that
\begin{itemize}
\item[1.] $U=\cup_{i}B_{i}$ with $B_{i}=
C_{1}\underline{B_{i}}$ and the balls $B_{i}$ have the bounded overlap property;
\item[2.] $r_{i}=r(B_{i})=\frac{1}{2}d(x_{i},F)$ and $x_{i}$ is 
the center of $B_{i}$;
\item[3.] each ball $\overline{B_{i}}=C_{2}\underline{B_{i}}$ intersects $F$ ($C_{2}=4C_{1}$ works).
\end{itemize}
Recall that the above balls are balls of $\Omega_{+}$, that is $\underline{B_{i}}= B(x_{i},r_{i}/C_{1}) \cap \Omega_{+}$, ${B_{i}}= B(x_{i},r_{i}) \cap \Omega_{+}$,  $\overline{B_{i}}= B(x_{i},r_{i}C_{2}) \cap \Omega_{+}$ and $x_{i}\in \Omega_{+}$ where $B(x,r)$ denotes an Euclidean open ball in $\RR^n$.  Condition (\ref{rb}) is nothing but the bounded overlap property of the $B_{i}$'s  and (\ref{eB}) follows from (\ref{rb}) and  (\ref{mO}). 

Since $ \lambda(\overline{B_{i}})
\leq C\lambda(B_{i})$ (the doubling property for $\Omega_{+}$) and  $\overline{B_{i}} \cap F
\neq \emptyset$ for all $i$, we have  
\begin{equation}\label{f}
\int_{B_{i}} (|f|+\frac{|f|}{r}+|\nabla f|)dx \leq
\int_{\overline{B_{i}}} (|f|+\frac{|f|}{r}+|\nabla f|) dx
\leq C \alpha\lambda(B_{i}).
\end{equation}

Let us derive some useful properties.
For $x\in U$, denote $I_{x}=\left\lbrace i:x\in B_{i}\right\rbrace$. By the bounded overlap property of the balls $B_{i}$, we have that $\sharp I_{x} \leq N$. Fixing $j\in I_{x}$ and using the properties of the $B_{i}$'s, we easily see that $\frac{1}{3}r_{i}\leq r_{j}\leq 3r_{i}$ for all $i\in I_{x}$. In particular, $B_{i}\subset 7B_{j}$ for all $i\in I_{x}$. We can deduce from that
  $|f_{B_{j}}-f_{B_{i}}|\leq Cr_{j} \alpha$ with $C$ independent of $i,j\in I_{x}$ and $x\in U$.
Indeed, we use that $B_{i}$ and $B_{j}$ are contained in $7 B_{j}$, Poincar\'e inequality $(P_{1})$ on balls of $\Omega_{+}$ (here $7B_{j}$), the comparability of $r_{i}$ and $r_{j}$, and the control of the gradient term in (\ref{f}). 

Let us now define the functions $b_{i}$ and  prove (\ref{eb}). Let $(\chi_{i})_{i}$ be  a partition of unity of $U$ subordinated to the covering $(\underline{B_{i}})$, such that for all $i$, $\chi_{i}$ is a Lipschitz function supported in $B_{i}$ with
$\displaystyle\|\,|\nabla \chi_{i}|\, \|_{\infty}\leq
\frac{C}{r_{i}}$. To this end it is enough to choose for $x\in \Omega_{+}$,  $\displaystyle\chi_{i}(x)=
\psi\Big(\frac{C_{1}d(x_{i},x)}{r_{i}}\Big)\Bigl(\sum_{k}\psi(\frac{C_{1}d(x_{k},x)}{r_{k}})\Bigr)^{-1}$, where $\psi$ is a smooth function, $\psi=1$ on $[0,1]$, $\psi=0$
on $[\frac{1+C_{1}}{2},+\infty[$ and $0\leq \psi\leq 1$. 	
We declare $B _{i}$ of type 1 if $4r_{i} \le d(B_{i},0)$ and of type 2 otherwise. Here $d(B_{i},0)$ is the distance from $B_{i}$ to 0. Indeed, it could well be that one $B_{i}$ even touches the origin. 
We set $b_{i}=(f-f_{B_{i}})\chi_{i}$ if $B_{i}$ is of type 1 and $b_{i}= f\chi_{i}$ if $f$ is of type 2. It is clear that $\supp b_{i} \subset B_{i}$.

 We first begin the proof of the estimates on $b_{i}$ by assuming $B_{i}$ of type 1. We remark that for all $x\in B_{i}$ we have $r_{i}\le |x|/4=r/4$ from the type 1. We have
\begin{equation*}
\int_{B_{i}} |b_{i}| dx
=\int_{B_{i}} |(f-f_{B_{i}})\chi_{i}| dx
\leq 2\int_{B_{i}}|f| dx
\leq C \alpha \lambda(B_{i}).
\end{equation*}
We applied  (\ref{f}) in the last step. For $\frac{|b_{i}|}{r}$, we use $1/r \leq 1/4r_{i}$ on $B_{i}$ and the Poincar\'e inequality $(P_{1})$ on $B_{i}$:
\begin{equation*}
\int_{B_{i}} \frac{|b_{i}|}{r} dx \leq
\int_{B_{i}} \frac{|f-f_{B_{i}}|}{4r_{i}} dx
\leq C \int_{B_{i}}|\nabla f| dx
\leq C \alpha \lambda(B_{i}),
\end{equation*}
using (\ref{f}) again.  For $\nabla b_{i}$, since $\nabla\bigl((f-f_{B_{i}})\chi_{i}\bigr)=\chi_{i}\nabla f
+(f-f_{B_{i}})\nabla\chi_{i}$, the Poincar\'e inequality $(P_{1})$ on $B_{i}$ and (\ref{f}) yield
\begin{align*}
\int_{B_{i}}|\nabla b_{i}|dx
 &\leq
C\left(\int_{B_{i}}|\chi_{i}\nabla f|dx
+\int_{B_{i}}|f-f_{B_{i}}|\,|\nabla \chi_{i}|dx \right)
\\
&\leq C \| \chi_{i}\|_{\infty}\alpha\lambda(B_{i})+ C\|\nabla \chi_{i}\|_{\infty}r_{i}\int_{B_{i}}|\nabla
f|dx
\\
&\leq C\alpha\lambda(B_{i}).
\end{align*}
Therefore (\ref{eb}) is proved if $B_{i}$ is a  type 1 ball.  

If $B_{i}$ is a type 2 ball, the control of $\int_{B_{i}} |b_{i}| dx$ and $\int_{B_{i}} \frac{|b_{i}|}{r} dx$ is direct from definition and (\ref{f}). For $\int_{B_{i}}|\nabla b_{i}|dx$, the only term requiring an argument is $\int_{B_{i}}|f\nabla \chi_{i}|dx$. We remark that the type 2 implies $|x| \le 2r_{i} + d(B_{i},0) \le 6r_{i}$ for $x\in B_{i}$. Hence, we have $|\nabla \chi_{i}(x)|\leq C/r_{i}\le  6C/r$ on $B_{i}$ and we can use (\ref{f}) again. Therefore (\ref{eb}) is proved if $B_{i}$ is a  type 2 ball.
Remark that we proved  
\begin{equation}
\label{bi}  \int_{B_{i}} {|b_{i}|}dx
\leq C\alpha {r_{i}}\lambda(B_{i})
\end{equation}
for all $i$ and also $|f_{B_{i}}| \le C\alpha {r_{i}}\lambda(B_{i})$ for type 2 balls.

Set $  g=f-\sum_{i}b_{i}$ so that (\ref{df}) is granted and it remains to establish (\ref{eg}). We begin by some observations. Since the sum is locally finite on $U$,  $g$ is defined  almost everywhere on $\Omega$ and $g=f$ on $F$. Observe that $g$ is a locally integrable function on $\Omega$. Indeed, let $\varphi\in L^{\infty}$ with compact support. Since $d(x,F)\geq r_{i}$ for $x \in \supp \,b_{i}$ and $\sum \lambda(B_{i}) \le C \lambda(U)$ by using doubling and the disjointness of the balls $\underline{B_{i}}$, we obtain, using (\ref{bi}),
\begin{align*} \int\sum_{i}|b_{i}|\,|\varphi|\,dx &\leq
\Bigl(\int\sum_{i}\frac{|b_{i}|}{r_{i}}\,dx\Bigr)\,\sup_{x\in
\Omega_+}\Bigl(d(x,F)|\varphi(x)|\Bigr)
\\
&
\leq
C\alpha\lambda(U) \sup_{x\in \Omega_+
}\Bigl(d(x,F)|\varphi(x)|\Bigr).
\end{align*}
Since $f\in L^{1}_{loc}(\Omega_{+})$, we deduce that $g\in L^{1}_{loc}(\Omega_{+})$\footnote{Note that since $b\in L^{1}$ in our case, we can say directly that $g\in L^{1}_{loc}$. However, this way of doing applies to the homogeneous case presented in Section 6.}.  We also note that by Lebesgue differentiation theorem, we have
\begin{equation}\label{Leb}
| f|  + \frac{|f|}{r}+|\nabla f| \le \alpha, \quad \lambda-\mathrm{a.e. \ on \ } F.
\end{equation}

We turn to proving the estimate  $\| g \|_{\infty} \le C\alpha$.  Using that $\sum_{i} \chi_{i}=1$ on $U$ and $0$ on $F$, we have $$ g=f\ind_{F}+\sum_{B_{i}\  \mathrm{type} \ 1}f_{B_{i}}\chi_{i}.$$

By (\ref{Leb}),  it remains  to estimate the series. By (\ref{f}), $ |f_{B_{i}}|\leq C \alpha$ and we conclude using $\sum_{i} \chi_{i} \le 1$.
 
 We continue with the estimate $\|{g}/{r}\|_{\infty} \le C\alpha$.  We use the same decomposition for $g$ as above. 
 On $F$, ${|g|}/{r}\leq \alpha$.  Let  $x,y\in B_{i}$ with $B_{i}$ a type 1 ball. We have $r_{i}\leq |y|/4$ by type 1 definition. We also have $|x-y|\leq 2r_{i}$. Hence we deduce that $|y|\le 2r_{i}+ |x| \le |y|/2+|x|$, so that $|y|\leq 2|x|$ and
  \begin{equation*}
 \frac{|f_{B_{i}}|}{|x|}\leq 2\barint_{\hspace{-5pt}B_{i}} \frac{|f(y)|}{|y|}dy\leq C \alpha.
 \end{equation*}
 Then let $x\in U$. Considering only the balls $B_{i}$ containing $x$, we have: 
 \begin{equation*}
 \frac{|g(x)|}{|x|} \leq \sum_{i\in I_{x}: B_{i}\  \mathrm{type} \ 1}\frac{|f_{B_{i}}|}{|x|}\chi_{i}(x)
 \leq  C \alpha \sum_{i\in I_{x}}\chi_{i}(x)= C\alpha.
 \end{equation*}

It remains to prove $\|\nabla g\|_{\infty}\leq C\alpha$. For that, we use the original representation of $g$, differentiate in the sense of distributions and calculate
\begin{equation*}
\nabla g = \nabla f -\sum_{i}\nabla b_{i}
=\nabla f-\bigg(\sum_{i}\chi_{i}\bigg)\nabla f -h
=\ind_{F}(\nabla f) -h
\end{equation*}
with
\begin{equation*}
h=\sum_{B_{i}\  \mathrm{type} \ 1}(f-f_{B_{i}})\nabla 
\chi_{i} + f \sum_{B_{i}\  \mathrm{type} \ 2}  \nabla 
\chi_{i}.
\end{equation*}
By (\ref{Leb}), $|\ind_{F}(\nabla f)| \leq \alpha$ $\lambda$-a.e.. We claim that a similar estimate holds for $h$, i.e. $|h(x)|\leq C\alpha$ for all $x\in \Omega_+$. For this, note first that $h$ vanishes on $F$.
Then fix  $x\in U$. Observe that 
$
\sum_{i}\nabla \chi_i(x)=0,
$
and by the definition of $I_{x}$, the sum reduces to $i\in I_{x}$.  Pick $j\in I_{x}$ with $B_{j}$ of type 2 if there is one such ball, otherwise any $j \in I_{x}$ will do. We have
$$
h(x)= \sum_{i\in I_{x } : B_{i}\  \mathrm{type} \ 1} (f_{B_{j}}-f_{B_{i}})\nabla \chi_i(x) + f_{B_{j}} \sum_{i\in I_{x } : B_{i}\  \mathrm{type} \ 2}\nabla 
\chi_{i}(x)
$$
because the difference with the previous equation is 
$$
(f(x) - f_{B_{j}}) \sum_{i \in I_{x}}\nabla \chi_i(x)=0.
$$
We have seen that  $|f_{B_{j}}-f_{B_{i}}|\leq Cr_{j} \alpha \leq 3Cr_{i}\alpha$ with $C$ independent of $i,j\in I_{x}$ and $x\in U$. 
 Since  $|\nabla \chi_i(x) | \leq C/r_{i}$ and  $I_{x}$ has cardinal bounded by $N$, we are done for the the first term in the right hand side. For the second term, either no $i\in I_{x}$ are such that $B_{i}$ is of type 2, in which case this term is 0. In the opposite case, we know that $|f_{B_{j}}|\leq C\alpha{r_{j}} $ since $B_{j}$ is of type 2 and we conclude using $|\nabla \chi_i(x) | \leq C/r_{i} \leq 3C/r_{j}$ for $i\in I_{x}$ and  that $I_{x}$ has cardinal bounded by $N$. From these estimates we deduce that $|\nabla g(x)|\leq C\alpha\; \lambda- a.e.$. 
 \qed
\\

\bigskip

We are now able to characterize the $K$-functional of interpolation between $\widetilde{H}_{1}^{1}(X)$ and $\widetilde{H}_{\infty}^{1}(X)$.
\begin{theorem} We have that
$$
K(f,t,\widetilde{H}_{1}^{1},\widetilde{H}_{\infty}^{1})
\sim t\left(f^{**}(t)+\left(\frac{|f|}{r}\right)^{**}(t)+|\nabla f|^{**}(t)\right)
$$
 for every $f\in \widetilde{H}_{1}^{1}(X)+\widetilde{H}_{\infty}^{1}(X)$ and $t>0$. The implicit constants are independent of $f$ and $t$.
\end{theorem}
\proof
The lower bound  follows from the fact that $K(g,t,L^1,L^\infty)\sim tg^{**}(t)$ for $g\in L^1+L^\infty$.
Now for the upper bound, consider first the case when $f\in \widetilde{H}_{p}^{1}(X)$. Identifying $f$ to its restriction to $\Omega$, write $f=f_+ + f_-$ and take the above Calder\'on-Zygmund decomposition for each $f_+$ and $f_{-}$ for $\alpha >0$ to be chosen. We obtain open subsets $U_\pm$ and functions $g_\pm, b_\pm$. We assume that $U_\pm$ are nonempty; the easy modifications otherwise are left to the reader. 

Here is the point of working with the $\widetilde{H}$ spaces instead of the $H$ spaces. As $g_{+}(0^+)=g_{-}(0^-)=0$, if we define $g=g_{+}$ on $\Omega_{+}$ and $g_{-}$ on $\Omega_{-}$, then $g$ can be extended to a Lipschitz function on $X=\overline{\Omega}$ with $\|\frac{g}{r}\|_{\infty}\leq C\alpha$. Hence $g\in \widetilde{H}_{\infty}^{1}(X)$ with norm controlled by $C\alpha$.  
 
 Therefore we can write $f\in \widetilde{H}_{p}^{1}(X)$ as $f=g+b$ with $b\in W_1^1(\Omega)=\widetilde{H}_{1}^{1}(X)$ and $g\in \widetilde{H}_{\infty}^{1}(X)$. We have $\|g\|_{\widetilde{H}_{\infty}^{1}(X)}\leq C \alpha$ and $\|b\|_{\widetilde{H}_{1}^{1}(X)}\leq C \alpha (\lambda(U_{+})+\lambda(U_{-}))$.  Let 
 $$
 \alpha_\pm(t)= \left(\mathcal{M}_{\Omega_\pm}(|f_\pm|+\frac{|f_\pm|}{r}+|\nabla f_\pm|)\right)^{*}(t), \quad \alpha= \max (\alpha_+(t), \alpha_-(t)).
 $$
 Remark that 
  \begin{align*}
  \alpha_+(t) 
  &\lesssim \left(|f_{+}|^{**}+\left(\frac{|f_{+}|}{r}\right)^{**}+|\nabla f_{+}|^{**}\right)(t)  \\
  &\lesssim \left(|f|^{**}+\left(\frac{|f|}{r}\right)^{**}+|\nabla f|^{**}\right)(t)
  \end{align*}
    where the implicit constant depends only on the doubling constant of $\Omega_+$.
We used the fact that
  $\left\lbrace x\in \Omega_{+};|f_{+}(x)|>\lambda\right\rbrace \subset \left\lbrace x\in \Omega;|f(x)|>\lambda\right\rbrace$, hence $f_{+}^{*}(t)\leq f^{*}(t)$. Similarly, $(\frac{f_+}{r})^{*}(t)\leq(\frac{f}{r})^{*}(t) $ and $|\nabla f_{+}|^{*}(t)\leq |\nabla f|^{*}(t)$.

 As $U_{+}$ is contained in $$ \left\lbrace x\in \Omega_+ ;  \mathcal{M}_{\Omega_+}(|f_+|+\frac{|f_+|}{r}+|\nabla f_+|)(x)>\alpha_+(t)\right\rbrace$$ we have $\lambda(U_+)\leq t$. Similarly we get $\lambda(U_{-})\leq t$. This yields
  \begin{align*}
K(f,t,\widetilde{H}_{1}^{1},\widetilde{H}_{\infty}^{1})&\leq \|b\|_{\widetilde{H}_{1}^{1}}+t\|g\|_{\widetilde{H}_{\infty}^{1}}
\\
&\leq C  t\left(f^{**}(t)+\left(\frac{|f|}{r}\right)^{**}(t)+|\nabla f|^{**}(t)\right).\label{DP4}
\end{align*}

For the general case when $f\in \widetilde{H}_{1}^{1}(X)+\widetilde{H}_{\infty}^{1}(X)$, we apply a similar argument to that of \cite{devore}  to obtain the upper bound. We  omit  details.
\qed
\proof[Proof of Theorem \ref{th:interpolationhptilde}] Set $
\widetilde{H}_{p,1}^{1}(X)=(\widetilde{H}_{1}^{1}(X),\widetilde{H}_{\infty}^{1}(X))_{1-1/p,p}.$
By the reiteration theorem,  it suffices to establish 
$\widetilde{H}_{p,1}^{1}(X)= \widetilde H^1_{p}(X)$ with equivalent norms.

First, from the Calder\'on-Zygmund decomposition, we have $\widetilde{H}_{p}^{1}(X)\subset
 \widetilde{H}_{1}^{1}(X)+\widetilde{H}_{\infty}^{1}(X)$ for $ 1<p<\infty$ where the inclusion is continuous.

From the previous results we have  that for $f\in \widetilde{H}_{1}^{1}(X)+\widetilde{H}_{\infty}^{1}(X)$
 \begin{align*}
\|f\|_{1-1/p,p}&\sim \left\lbrace\int_{0}^{\infty}\left(|f|^{**}(t)+\left(\frac{|f|}{r}\right)^{**}+|\nabla f|^{**}(t)\right)^{p}dt\right\rbrace^{1/p}
\\
&\sim \|f^{**}\|_{p}+\left\|\left(\frac{|f|}{r}\right)^{**}\right\|_{p}+
+\| \, |\nabla f|^{**} \|_{p}
\\
&\sim \|f\|_{p}+\left\|\frac{f}{r}\right\|_{p}
+\| \,|\nabla f|\,\|_{p}
\\
&\sim \| f\|_{\tilde{H}_{p}^{1}},
\end{align*}
where we used that for $l>1$, $\|f^{**}\|_{l}\sim \|f\|_{l}$.
\qed

\section{Restriction/Extension from/to $\RR^n$}
We study the restriction operator onto $\Omega$ and 
 construct an extension that is $p$ independent.  The subject of restriction and extension has been very studied for Sobolev spaces on domains (that is, connected open sets). For some definite answers see \cite{hajlasz3} and the references therein. But recall that the double cone is not a domain. For references closer to what we are doing here, see \cite{R} which considers the Bessel-Sobolev spaces in subsets of $\RR^n$ and \cite{S} which treats the Hajlasz-Sobolev spaces in spaces of homogeneous type. Let us first state our result.

\begin{theorem}   Let $1\le p \le \infty$.
\begin{itemize}
\item
The restriction operator is bounded from $W^1_{p}(\RR^n)$ into $H^1_{p}(X)$. Further, it is onto for $p\ne n$ and for $p=n$, its range is  $\widehat{H}_n^1(X)$.
\item There  exists a linear extension operator  $E$ that is bounded from $H^1_{p}(X)$ to 
$W^1_{p}(\RR^n)$ if $p\ne n$ and from $\widehat{H}_n^1(X)$  to 
$W^1_{n}(\RR^n)$.
\end{itemize}
\end{theorem}

The interesting part of this result is $p=n$. Observe also that this shows that $H^1_{n}(X)$ does not have the extension property.

In accordance with the cited references, the second item is closely related to ontoness in the first. It could well be a direct consequence but we prefer producing  an explicit extension operator. 

In \cite{S}, the author studies restrictions onto regular sets in a space of homogeneous type. Note that the double cone is a regular set in $\RR^n$. For such sets, he proves the following interesting theorem (Theorem 1.3). For all $1<p\le \infty$, the restriction of $M^{1,p}(\RR^n)$ to $\Omega$ equals  $M^{1,p}(X)$ and there exists a linear continuous extension from the latter space to the first one. Here, $M^{1,p}$ is the Hajlasz-Sobolev space. Given the fact that $W^1_{p}(\RR^n)=M^{1,p}(\RR^n)$ with equivalent norms for $1<p\le \infty$  (\cite{hajlasz2}), this  implies that $M^{1,p}(X)$ interpolate for $1<p\le \infty$. But combining this with our result gives the following corollary.

\begin{corollary} For $1<p\le \infty$, $M^{1,p}(X)=H^1_{p}(X)$ if $p\ne n$ and $M^{1,n}(X)= \widehat{H}_n^1(X)$.
\end{corollary}

The result for $p>n$ is known (see \cite{hajlasz1}): it follows from the fact that $X$ is complete and satisfies Poincar\'e inequality for any $p>n$. Not much more can be said in general  without these two conditions so this corollary is seemingly new  for $p\le n$.  The interesting case is the identification for $p=n$: $M^{1,n}(X)$ is a strict subspace of $H^{1}_{n}(X)$.

This study of restriction/extension properties  uses the decomposition into radial and anti-radial parts defined earlier. It also possible  to reprove the interpolation property of the $H^1_p(X)$ spaces by this method. 

We first study the restriction operator, then construct the extension.  We  next prove the ontoness and conclude with the application to interpolation.

\subsection{Restriction}
The restriction operator $R$ is defined by $R(f)=f|_\Omega$. Let $1\leq p< \infty$. It is obvious that if $f\in W_p^1(\mathbb{R}^n)$ then $R(f) \in W^1_p(\Omega)$, and that  $R: W_p^{1}(\mathbb{R}^n)\rightarrow W_p^1(\Omega)$ is bounded.  As $C_0^\infty(\RR^n)$ is dense in $W_p^1(\RR^n)$, the range is contained in $\widetilde W^1_p(\Omega)=H^1_p(X)$. 

For $p= n$, we show that $R$ maps into $\widehat H^1_{n}(X)$. Let $f\in W^1_{n}(\RR^n)$ and let $g=R(f)$. Since we already know that $g\in H^1_{n}(X)$, it remains  to show that $g_{a}/r\in L^n(X)$. Write  $f=f_\rho+f_\alpha$ where  $f_\rho(x)$ is here the average of $f$ on the whole sphere of radius $|x|$. Identifying $f_{\rho}$ with its restriction to $\Omega$, we see that $(f_{\rho})_{r}=f_{\rho}$ and $(f_{\rho})_{a}=0$. Thus, if we write $g=g_{r}+g_{a}$, we conclude that $g_{a}=( {f_{\alpha}|_{\Omega}})_{a}$, i.e $g_{a}(x)= f_{\alpha}(x) - \barint_{\Omega\cap S_{|x|}}f_{\alpha} \, d\sigma(\theta)$ for $x\in \Omega$.
Thus
\begin{align*}
 \int_{\Omega\cap S_{|x|}} |g_{a}|^n \, d\sigma(\theta)& \leq 2^n  \int_{\Omega\cap S_{|x|}} |f_{\alpha}|^n \, d\sigma(\theta)  \\
 &
 \leq  2^n \int_{S_{|x|}} |f_{\alpha}|^n \, d\sigma(\theta) 
 \\
 &\leq C |x|^n 
 \int_{S_{|x|}}|\nabla_{\theta}f_\alpha|^{n}d\sigma(\theta)
 \end{align*}
 where the last inequality is Poincar\'e inequality $(P_{n})$ on the sphere $S_{|x|}$ since $f_{\alpha}$ has mean value 0 on it, and $\nabla_{\theta}$ is the tangential gradient (see the last remark in Section 3). Since $r(x)=|x|$ it follows that  
\begin{align*}
\int_{\Omega}\left|\frac{g_a}{r}\right|^{n}dx&=\int_{0}^{+\infty}\int_{\Omega\cap S_{r}}\left|\frac{g_a}{r}\right|^{n}d\sigma(\theta)dr
\\
&\leq C\int_{0}^{+\infty} \int_{S_{r}}|\nabla_{\theta}f_\alpha|^{n}d\sigma(\theta)dr
\\
&= C\int_{0}^{+\infty} \int_{S_{r}}|\nabla f_{\alpha}|^{n}d\sigma(\theta)dr
\\
&\leq C\int_{\RR^n}|\nabla f|^{n}dx<\infty.
\end{align*} 
We used that $\nabla_{\theta}f_{\alpha}=\nabla_{\theta} f$ since $f_{\rho}$ does not depend on $\theta$ and then $|\nabla f_{\alpha}| \le |\nabla f|$. Therefore  $R(f)=g \in  \widehat{H}_n^1(X)$.  

For $p=\infty$, it is obvious that $R$ is bounded from $W^1_\infty(\RR^n)$ into $H^1_\infty(X)$ (and it is onto by Whitney's extension theorem).

\subsection{Extension}
Let $1\leq p \leq \infty$. Let $f\in H_{p}^1(X)$ if $p\ne n$ (resp. $f\in \widehat{H}_n^1(X)$). Write $f=f_r+f_a$ as in Section 3. Lemma \ref{lemma:ra} yields $f_r \in W_p^1(\mathbb{R}^n)$ and $f_a\in \widetilde{H}_p^1(X)$. It remains to extend $f_a$. We write $f_a=f_a|_{\Omega_+}+f_a|_{\Omega_{-}}=f_{a+}+f_{a-}$. We treat $f_{a+}$, the same analysis applying to  $f_{a-}$. Our strategy is to enlarge $\Omega_{+}$ to a slightly bigger half-cone $\widetilde{\Omega}_+$, then we map them with a smooth bilipschitz map onto the upper half-space $\mathbb{R}^{n}_+$ and a slightly bigger version $\widetilde{\mathbb{R}^{n}_{+}}$. We use the reflection principle for Sobolev spaces to extend from the upper-half space to the full space, localise on $\widetilde{\mathbb{R}^{n}_{+}}$  with a homogeneous cut-off and finish by mapping back onto $\widetilde{\Omega}_+$. There are many ways to do that. Here are the details. 

Let $w<\pi/2$ be the half-angle of $\Omega_{+}$ with respect to the vertical axis and  let  $\widetilde{\Omega}_+$ be an open half-cone with half-angle $\omega +\varepsilon<\pi/2$ for some small $\varepsilon$ with same vertex point and rotation axis. Using the spherical angle $\theta \in [0,\pi)$ defined by $ \theta= \arccos\frac{x_{n}}{|x|}$ by writing $x=(x',x_{n}), x'\in \RR^{n-1}, x_{n}\in \RR$, the map  $\psi_+(x)=y$ with $y=\big( \frac{\sin (2\omega \theta/\pi)}{\sin \theta} \, x', \frac{\cos (2\omega \theta/\pi)}{\cos \theta} \, x_{n}\big)$ is a smooth bilipshitz  (this is somewhat tedious to check carefully but it reduces to a planar estimate) map, leaving the norm invariant ($|\psi_{+}(x)|=|x|$), from $\mathbb{R}^{n}_+$ onto $\Omega_+$ and from $\widetilde{\mathbb{R}^{n}_{+}}$ onto $\widetilde{\Omega}_+$  where $\widetilde{\mathbb{R}^{n}_{+}}$ is a half-cone with half-angle $\frac{\pi (w+\varepsilon)}{2\omega} >\frac{\pi}{2}$. We consider now the even extension $\zeta_+: W_{p}^{1}(\mathbb{R}^{n}_+)\rightarrow W_p^1(\mathbb{R}^{n})$. Let $m_+\in C^{\infty}(\mathbb{R}^{n}\setminus \{0\})\cap L^\infty(\RR^n)$, homogeneous of degree 0,  such that  $m_+=1$ on $\mathbb{R}^{n}_+$ and $\supp m_+ \subset \widetilde{\mathbb{R}^n_{+}}\cup\{0\}$. With these ingredients we define the extension $\xi_+(f_{a+})$ of $f_{a+}$ as
$$\xi_+(f_{a+})=\big[m_+\zeta_+(f_{a+}\circ \psi_+)] \circ \psi_+^{-1}.$$

It readily follows from the properties of $m_+$ that 
$$
\left\| \frac{m_+ \,g}{r}\right\|_p \lesssim \left\| \frac{g}{r}\right\|_p
$$
and
$$
\|\, |\nabla (m_+g)| \|_p \lesssim   \|\, |\nabla g|\|_p +\left\| \frac{g}{r}\right\|_p
$$
for all $g\in W^1_p(\RR^n)$. Using this fact,  that bilipschitz maps preserve Sobolev spaces and density of Lipschitz functions, we obtain  
 that $\xi_+(f_{a+})\in W_p^1(\mathbb{R}^n)$ with $\supp \xi_+(f_{a+})\subset \widetilde{\Omega}_{+} \cup \{0\}$ and $\|\xi_+(f_{a+})\|_{W_p^1(\mathbb{R}^n)}\leq C\|f_{a+}\|_{\widetilde{H}_{p}^1(X)}$. 
We conclude that $\xi(f_a)=\xi_+(f_{a+})+\xi_-(f_{a-})\in W_p^1(\mathbb{R}^n)$ is an extension of $f_a$ to $W^1_p(\RR^n)$. 

Therefore, $E$ defined by $$E(f)=f_r+\xi(f_a)$$ is an extension of $f$ to $W^1_p(\RR^n)$. We have shown that  the map $E$ is $ H_p^1(X)\rightarrow W_p^1(\RR^{n})$-bounded if $p\neq n$  and $\widehat{H}_n^1(X)\rightarrow W_n^1(\mathbb{R}^n) $-bounded if $p=n$.

\subsection{Ontoness and relation to  interpolation for $H_p^1(X)$}
From the previous subsections, we deduce that $R\circ E$ operates boundedly on $H_{p}^{1}(X)$ for $p\neq n$ and on  $\widehat{H}_n^1(X)$  as the identity map. In particular, $R$ acting on $W^1_{p}(\RR^n)$ is onto $H_{p}^{1}(X)$ for $p\neq n$ and onto  $\widehat{H}_n^1(X)$ for $p=n$. Using the preservation of interpolation properties for retract diagrams, it follows that  
$$(H_1^1(X),H_{\infty}^{1}(X))_{1-1/p,p}=R(W_{1}^1(\mathbb{R}^{n}),W_{\infty}^1(\mathbb{R}^{n}))_{1-1/p,p})=R(W_p^1(\mathbb{R}^n)).$$
Therefore $(H_1^1(X),H_{\infty}^{1}(X))_{1-1/p,p}=H_p^1(X)$ for $p\neq n$ and $\widehat{H}_n^1(X)$ for $p=n$.

\section{Homogeneous versions} 

Homogeneous Sobolev spaces are defined up to a constant, removing control on the $L^p$ norms on $f$. Since the vertex point plays a specific role, it is best here to fix the floatting constant by imposing control at this vertex point. We adopt the following definitions.
Let  $\Lip^{0}(X)$ be the space of Lipschitz functions in $X$ vanishing at 0. For $1 \leq p \le \infty$,  we set 
$$E_p =\{ f \in \Lip^{0}(X) \, ; \, \|\, |\nabla f|\|_{L^p(X)} <\infty\},
$$
$$
\widetilde E_p= \{ f \in \Lip^{0}(X) \, ; \, \|\, |\nabla f|\|_{L^p(X)} + \|f/r\|_{L^p(X)} <\infty\}.
$$
Then $E_p$ and $\widetilde E_p$ are normed spaces and we call $ \mH^1_p(X)$ and $\widetilde \mH^1_p(X)$  their completions. Clearly  $\widetilde E_\infty =\widetilde\mH^1_\infty(X)= E_\infty=\mH^1_\infty(X)=  \Lip^{0}(X)$.

It is easy  to show that $ \mH^1_p(X)$ is composed of locally $p$-integrable functions. 
For $p>n$, one has  \eqref{eq:Morrey} from the Morrey embedding and $f(0)=0$.

It is clear that $\widetilde \mH^1_p(X) \subset \mH^1_p(X)$ but  for $1\le p\le n$ the inclusion is strict\footnote{In contrast with the inhomogeneous case for $p<n$.}. Indeed a Lipschitz function supported away  from 0 which agrees with $e^{i|x|^{-\alpha}}$ for $|x|\ge 1$  satisfies $\|f/r\|_p=\infty$ and  belongs to $E_p$ if $\alpha>0$ is large enough. For $p>n$, the inclusion is an equality as we shall see.

\begin{lemma} For $1\le p <\infty$, $\Lip^0(X)\cap \Lip_0(X)$ is dense in $\widetilde \mH^1_p(X)$.
\end{lemma}

\proof If $f \in \widetilde E_p$, consider $f_k=f \chi(r/k)$, $k \in \NN^*$, where $\chi: [0,\infty) \to [0, 1]$ is a smooth function which is 1 on $[0,1]$ with support in $[0,2]$. It is easy to show that $\| \, |\nabla (f-f_k)|\|_p$ and $\|(f-f_k)/r\|_p$ tend to 0 as $k$ tends to $\infty$.\qed

\begin{remark} From there, one can see that the restrictions to $\Omega$ of  functions in $C_0^\infty(\RR^n)$ that vanish at $0$ form a dense subspace of $\widetilde \mH^1_p(X)$.
\end{remark}

\begin{corollary} \begin{itemize}
\item For $1\le p < n$, \eqref{uc} holds on $\widetilde \mH^1_p(X)$.
\item For $n< p <\infty$,   \eqref{uc} holds on $ \mH^1_p(X)$  and  $ \mH^1_p(X)=\widetilde \mH^1_p(X)$.
\end{itemize}
\end{corollary}

\proof  Assume  first that $1\le p<n$.  Then by the previous lemma, one can assume that $f\in \Lip_{0}(X)$ for which the argument of \eqref{uc} applies. 

Assume now $n<p<\infty$. Let $f \in \Lip^0(X)$. For $0<\epsilon <R <\infty$, set
$A=\int_{\Omega_{+}\cap\{R>|x|>\epsilon\}}\big|\frac{f}{r}\big|^{p}dx$. Then argue as in the proof  of \eqref{uc}. In the integration by parts, one picks an extra term  which has a negative sign because $n-p<0$. Thus one can cancel it and obtain $A \le C \| | \nabla f|\|_{L^p(X_+)}^p$ with  $C$ independent of $\epsilon, R$. Taking limits and doing the same thing on $\Omega_-$ shows that $f\in \widetilde E_p$ and we are done. 
\qed

The first item also show that the closure in $\mH^1_{p}(X)$ of $\Lip^0(X)\cap \Lip_0(X)$ is 
$\widetilde \mH^1_p(X)$.

\begin{theorem} The family $(\widetilde \mH^1_p(X))_{1\le p \le \infty}$ is an interpolation family for the real method. Hence the same is true for $( \mH^1_p(X))_{n< p \le \infty}$.
\end{theorem}

The proof for the spaces $\widetilde \mH^1_p(X)$ is a minor adapatation of the one of Theorem \ref{th:interpolationhptilde} and is left to the reader. The second point follows from the above corollary.

Interpolation  for the spaces $\mH^1_p(X)$ for $p\le n$ is unclear.

\section{Some remarks and generalizations}

\begin{remark} (some explicit extensions) There are many extension operators. The following example was communicated to us by Michel Pierre. For the (double) cone of $\RR^2$ consisting of the 2 quadrants defined by $xy>0$, then one can take  
  $$
Ef(x,y)=\begin{cases} f(x,y), & \text{if} \ xy > 0, \\
\dfrac{x^2   f(x,-y) + y^2  f(-x,y) }{x^2+y^2}, & \text{if} \ xy < 0.
\end{cases}
 $$

\end{remark}

\begin{remark} (central role of the vertex) The analysis in this article  does not use the fact that the cone $\Omega$ is symmetric under $x\mapsto -x$ and also does not use the specific opening angle. This means that the upper and lower (open) half-cones can be replaced by two half-cones located independently of one another provided they share the same vertex and that they are strictly separated by a hyperplane passing through the vertex and not containing any direction of the boundaries. Also the (finite) number of  disjoint half-cones is not limited to 2 provided each pair satisfies the above requirements. 
\end{remark}

\begin{remark} (other types of cones) The half-cones can be replaced by $\RR_{+}^*\times N$ where $N$ is a Lipschitz domain on the unit sphere. On such domains, one has Poincar\'e inequalities with any exponents (adapt the proof sketched in the last remark of Section 3) and this allows to adapt the arguments. 
\end{remark}

\begin{remark} (local geometry) Of course, the analysis done with inhomogeneous norms is stable  by  (smooth) truncation of the cone away from the vertex point. For example, if one wants to work on  a truncated cone by requiring $r<1$, then one can use local variants as in  Badr's thesis \cite{badr}. Details are left to the reader.
\end{remark}

\bibliographystyle{amsplain}
\addcontentsline{toc}{section}{
\end{document}